\documentclass{elsart}

\usepackage{amssymb,amsfonts,amsxtra}
\usepackage{dsfont}
\usepackage[all]{xypic}
\xyoption{dvips}

\newcommand{\dd}{\mathrm{d}}
\newcommand{\df}{\dd{f}}
\newcommand{\dt}{\partial_{t}}
\newcommand{\dta}{\partial_{\tau}}
\newcommand{\dth}{\partial_{\theta}}

\newcommand{\ee}{\mathrm{e}}
\newcommand{\ii}{\mathrm{i}}
\newcommand{\inv}[1]{#1^{-1}}
\newcommand{\ol}[1]{{\overline{#1}}}
\newcommand{\pd}[1][]{\partial_{#1}}

\newcommand{\sO}{\mathcal{O}}
\newcommand{\ul}[1]{{\underline{#1}}}
\newcommand{\upd}[1][]{\ul\partial_{#1}}
\newcommand{\wh}{\widehat}
\newcommand{\xymat}{\SelectTips{cm}{}\xymatrix}

\newcommand{\C}{\mathds{C}}
\newcommand{\Cmth}{\C\{\!\{\theta\}\!\}}

\newcommand{\Cta}{\C[\tau]}
\newcommand{\Ctad}{\Cta\langle\dta\rangle}
\newcommand{\Ct}{\C[t]}
\newcommand{\Ctd}{\Ct\langle\dt\rangle}
\newcommand{\Cth}{\C[\theta]}
\newcommand{\Cthd}{\C[\theta]\langle\dth\rangle}
\newcommand{\Ctht}{\C[\theta]\langle t\rangle}
\newcommand{\Ctt}{\C[\tau,\theta]}
\newcommand{\Cttt}{\C[\tau,\theta]\langle t\rangle}
\newcommand{\CX}{\C[\ul{x}]}
\newcommand{\CXth}{\C[\ul x,\theta]}
\newcommand{\CXtt}{\C[\ul x,\tau,\theta]}

\newcommand{\N}{\mathds{N}}

\newcommand{\Q}{\mathds{Q}}
\newcommand{\Z}{\mathds{Z}}

\DeclareMathOperator{\crt}{C}
\DeclareMathOperator{\dsc}{D}
\DeclareMathOperator{\gr}{gr}
\DeclareMathOperator{\lead}{lead}
\DeclareMathOperator{\lexp}{lexp}
\DeclareMathOperator{\rk}{rk}
\DeclareMathOperator{\spec}{spec}
\DeclareMathOperator{\End}{End}
\DeclareMathOperator{\GL}{GL}
\DeclareMathOperator{\HH}{H}

\DeclareMathOperator{\NF}{NF}
\DeclareMathOperator{\RR}{R}
\DeclareMathOperator{\GB}{GB}

\begin{document}

\begin{frontmatter}

\title{Good bases for tame polynomials}
\author{Mathias Schulze\thanksref{a}}
\address{Department of Mathematics\\University of Kaiserslautern\\67663 Kaiserslautern\\Germany}
\ead{mschulze@mathematik.uni-kl.de}
\thanks[a]{The author is grateful to Claude Sabbah for drawing his attention to the subject and to Gert--Martin Greuel for valuable hints and discussions.}

\begin{abstract}
An algorithm to compute a good basis of the Brieskorn lattice of a cohomologically tame polynomial is described.
This algorithm is based on the results of C.~Sabbah and generalizes the algorithm by A.~Douai for convenient Newton non--degenerate polynomials.
\end{abstract}

\begin{keyword}
tame polynomial \sep Gauss--Manin system \sep Brieskorn lattice \sep V--filtration \sep mixed Hodge structure \sep monodromy \sep good basis
\MSC 13N10 \sep 13P10 \sep 32S35 \sep 32S40
\end{keyword}

\end{frontmatter}

\section*{Introduction}

Let $\xymat{f:\C^{n+1}\ar[r]&\C}$ with $n\ge1$ be a cohomologically tame polynomial function \cite{Sab98b}.
This means that no modification of the topology of the fibres of $f$ comes from infinity.
In particular, the set of critical points $\crt(f)$ of $f$ is finite.
Then the reduced cohomology of the fibre $\inv f(t)$ for $t\notin\crt(f)$ is concentrated in dimension $n$ and equals $\C^\mu$ where $\mu$ is the Milnor number of $f$.
Moreover, the $n$--th cohomology of the fibres of $f$ forms a local system $H^n$ on $\C\backslash\dsc(f)$ where $\dsc(f)=f(\crt(f))$ is the discriminant of $f$.
Hence, there is a monodromy action of the fundamental group $\Pi_1(\C\backslash\dsc(f),t)$ on $H^n_t$.

The Gauss--Manin system $M$ of $f$ is a regular holonomic module over the Weyl algebra $\Ctd$ with associated local system $H^n$ on $\C\backslash\dsc(f)$.
The Fourier transform $G:=\wh M$ of $M$ is the $\Ctad$--module defined by $\tau:=\dt$ and $\dta=-t$.
The monodromy $T_\infty$ of $M$ around $\dsc(f)$ can be identified with the inverse of the monodromy $\wh T_0$ of $G$ at $0$.
It turns out that $\dt$ is invertible on $M$ and hence $G$ is a $\Ctt$--module where $\theta:=\inv\tau$.
A finite $\Cta$-- resp. $\Cth$--submodule $L\subset G$ such that $L[\theta]=G$ resp. $L[\tau]=G$ is called a $\Cta$-- resp. $\Cth$--lattice.
The regularity of $M$ at $\infty$ implies that $G$ is singular at most in $\{0,\infty\}$ and where $0:=\{\tau=0\}$ is regular and $\infty:=\{\theta=0\}$ of type $1$.
In particular, the V--filtration $V_\bullet$ on $G$ at $0$ consists of $\Cta$--lattices.

The Brieskorn lattice $G_0\subset G$ is a $t$--invariant $\Cth$--submodule of $G$ such that $G=G_0[\tau]$.
C.~Sabbah \cite{Sab98b} proved that $G_0$ is a free $\Ct$-- and $\Cth$--module of rank $\mu$.
In particular, $G$ is a free $\Ctt$--module of rank $\mu$.
By definition, the spectrum of a $\Cth$--lattice $L\subset G$ is the spectrum of the induced V--filtration $V_\bullet(L/\theta L)$ and the spectrum of $f$ is the spectrum of $G_0$.

C.~Sabbah \cite{Sab98b} showed that there is a natural mixed Hodge structure on the moderate nearby cycles of $G$ with Hodge filtration induced by $G_0$.
This leads to the existence of good bases of the Brieskorn lattice.
For a basis $\ul\phi=\phi_1,\dots,\phi_\mu$ of a $t$--invariant $\Cth$--lattice,
\[
t\circ\ul\phi=\ul\phi\circ(A^{\ul\phi}+\theta^2\dth)
\]
where $A^{\ul\phi}\in\Cth^{\mu\times\mu}$.
A $\Cth$--basis $\ul\phi$ of $G_0$ is called good if $A^{\ul\phi}=A^{\ul\phi}_0+\theta A^{\ul\phi}_1$ where $A^{\ul\phi}_0,A^{\ul\phi}_1\in\C^{\mu\times\mu}$,
\[
A^{\ul\phi}_1=\begin{pmatrix}
\alpha_1 & & \\
 & \ddots & \\
 & & \alpha_\mu
\end{pmatrix}
\]
and $\phi_i\in V_{\alpha_i}G_0$ for all $i\in[1,\mu]$.
One can read off the monodromy $T_\infty=\inv{\wh T_0}$ from $A^\ul\phi$ immediately.
The diagonal $\ul\alpha=\alpha_1,\dots,\alpha_\mu$ is the spectrum of $f$ and determines with $\gr^V_1A_0$ the spectral pairs of $f$.
The latter correspond to the Hodge numbers of the above mixed Hodge structure.

Analogous results to those above were first obtained in a local situation where $\xymat{f:(\C^n,\ul 0)\ar[r]&(\C,0)}$ is a holomorphic function germ with an isolated critical point \cite{Bri70,Seb70,Ste76,Pha79,Var82a,Sai89}.
In this situation, the role of the Fourier transform is played by microlocalization and the algorithms in \cite{Sch02c,Sch04a} compute $A_0$ and $A_1$ for a good $\Cmth$--basis of the (local) Brieskorn lattice.
But \cite{Sch04a} and \cite[7.4--5]{Sch02c} do not apply to the global situation.

A.~Douai \cite{Dou99} explained how to compute a good basis of $G_0$ if $f$ is convenient and Newton non--degenerate using the equality of the V-- and Newton filtration \cite{KV85,Sab98b} and a division algorithm with respect to the Newton filtration \cite{Dou93,BGM89}.

The intention of this article is to describe an explicit algorithm to compute a good basis of $G_0$ for an arbitrary cohomologically tame polynomial $f$.
This algorithm is based on the following idea:

Let $\ul x=x_0,\dots,x_n$ be a coordinate system on $\C^{n+1}$.
Then the Brieskorn lattice $G_0$ can be identified with the quotient
\[
\CXth\Big/\sum_{i=0}^n(\pd[x_i](f)-\theta\pd[x_i])(\CXth)
\]
of non--finite $\Cth$--modules.
The degree with respect to $\ul x$ defines an increasing filtration $\CXth_\bullet$ by finite $\Cth$--modules on $\CXth$ and hence 
\[
G_0^{k,l}:=\CXth_k\Big/\Bigl(\CXth_k\cap\sum_{i=0}^n(\pd[x_i](f)-\theta\pd[x_i])(\CXth_l)\Bigr)
\]
are finite $\Cth$--modules.
For $k\gg0$ and $l\gg0$, $G_0^{k,l}=G_0$ by the finiteness of $G_0$.
But, a priori, there is no bound for these indices.

By Gr\"obner basis methods, one can compute cyclic generators $\ul\phi$ of a $t$--invariant $\Cth$--sublattice $G_0^{k,l}\subset G_0$.
By an argument of A.~Khovanskii and A.~Varchenko \cite{KV85}, $G_0^{k,l}=G_0$ if and only if the mean values of the spectra coincide.
By the $t$--invariance of $G_0^{k,l}$, one can compute the spectrum of $G_0^{k,l}$ like that of $G_0$ below.
The mean value of the spectrum of $G_0$ is known to be $\frac{n+1}{2}$.
So if the mean value of the spectrum of $G_0^{k,l}$ is not $\frac{n+1}{2}$ then one has to increase $k$.
This process terminates with $G_0^{k,l}=G_0$.

Then one can compute $A^{\ul\phi}$ for the $\Cth$--basis $\ul\phi$ of $G_0$.
By a saturation process, one can compute the V--filtration and, by a Gr\"obner basis computation, the spectrum of $G_0$ and the Hodge filtration.
Then one can compute a $\Ctt$--basis of $G$ which is compatible with the V--filtration refined by an opposite Hodge filtration.
In terms of this basis, one can compute a good basis of $G_0$ by a simultaneous normal form computation and basis transformation.

We denote rows vectors $\ul v$ by a lower bar and column vectors $\ol v$ by an upper bar.
In general, lower indices are column indices and upper indices are row indices.
We denote by $\{M\}$ the set and by $\langle M\rangle R$ the $R$--linear span of the of columns of a matrix $M$.
We denote by $\lead$ the leading term and by $\lexp$ the leading exponent with respect to a monomial ordering.
We denote by $E$ the unit matrix and by $\ol e_i$ the $i$th unit vector.

\section{Gauss--Manin system}

Let $\xymat{f:\C^{n+1}\ar[r]&\C}$ with $n\ge1$ be a polynomial function.
Let $\sO$ be the sheaf of regular functions and $(\Omega^\bullet,\dd)$ the complex of polynomial differential forms on $\C^{n+1}$.
Then the Gauss--Manin System $f_+\sO$ of $f$ is represented by the complex of left $\Ctd$--modules
\[
(\Omega^{\bullet+n+1}[\dt],\dd-\dt\dd f)
\]
\cite[15]{Pha79} and has regular holonomic cohomology \cite[VII.12.2]{Bor87}.
The coefficients of the differentials are the differentials of the complexes $(\Omega^\bullet,\dd)$ and $(\Omega^\bullet,\dd f)$.

\begin{lem}[Poincar\'e Lemma]\label{38}
The complex of $\C$--vector spaces
\[
\xymat{0\ar[r]&\C\ar[r]&(\Omega^\bullet,\dd)\ar[r]&0}
\]
is exact \cite[Ex.~16.15]{Eis96}.
\end{lem}

From now on, we assume that set of critical points $\crt(f)$ of $f$ is finite.
Then the following lemma holds.

\begin{lem}[De~Rham Lemma]\label{39}\
\begin{enumerate}
\item $\HH^k(\Omega^\bullet,\dd f)=0$ for $k\ne n+1$.
\item $\dim_\C\HH^{n+1}(\Omega^\bullet,\dd f)<\infty$.
\end{enumerate}
\end{lem}
\begin{pf}
If $\crt(f)$ is finite then
\[
\CX/\langle\upd(f)\rangle\cong\Omega^{n+1}/\dd f\wedge\Omega^{n-1}=\HH^{n+1}(\Omega^\bullet,\dd f)
\]
is a finite $\C$-vector space and hence $\upd(f)$ is a regular sequence in $\CX$.
Then the cohomology of the Koszul complex $(\Omega^\bullet,\dd f)$ is concentrated in dimension $n+1$ \cite[Cor.~17.5]{Eis96}. 
\end{pf}

The image $M_0$ of $\Omega^{n+1}$ in $M:=\HH^0(f_+\sO)$ is the key for an algorithmic approach to the Gauss--Manin system.
It determines the differential structure of $M$ and it can be identified with a quotient of $\CX$.

\begin{prop}\label{19}\
\begin{enumerate}
\item $\HH^k(f_+\sO)=0$ for $k\notin\{-n,0\}$.
\item $\dt$ is invertible on $M$.
\item $M_0=\Omega^{n+1}/\df\wedge\dd\Omega^{n-1}$.
\end{enumerate}
\end{prop}
\begin{pf}
This follows from Lemma \ref{38} and \ref{39} and \cite[15.2.2]{Pha79}.
\end{pf}

The Fourier transform $\wh M$ of $M$ is the left $\Ctad$--module defined by the isomorphism
\[
\tau:=\dt,\quad\dta:=-t
\]
of $\Ctad$ and $\Ctd$ \cite[2.1]{Sab98a}.
By Proposition \ref{19}, $M$ is a $\Cthd$--module with
\[
\theta:=\inv\tau,\quad\dth:=-\tau^2\dta
\]
and $M_0$ is a $\Cth$-submodule.
Note that $t=\theta^2\dth$.
\begin{defn}
Let $G$ be the $\Cthd$-module $\wh M$.
Then the Brieskorn lattice $G_0$ of $f$ is the $\Ctht$-submodule $M_0$ of $G$.
\end{defn}

Since $M$ is regular at $\infty$, $G$ is singular at most in $\{0,\infty\}$ where $0:=\{\tau=0\}$ is regular and $\infty:=\{\theta=0\}$ of type $1$ \cite[V.2.a]{Sab02}.

From now on, we assume that $f$ is cohomologically tame.
By definition \cite[8]{Sab98b}, this means that there is a compactification
\[
\xymat{
\C^{n+1}\ar@{^(->}[r]^-j\ar[rd]^-f&\ol{\C^{n+1}}\ar[d]^-{\ol f}\\
&\C
}
\]
where $\ol{\C^{n+1}}$ is quasi--projective and $\ol f$ is proper such that, for all $t\in\C$, the support of the vanishing cycle complex $\phi_{\ol f-t}\RR j_*\Q$ is a finite subset of $\C^{n+1}$.
In particular, $\crt(f)$ is finite and hence, by Lemma \ref{39}, the Milnor number
\[
\mu:=\dim_\C\HH^{n+1}(\Omega^\bullet,\dd f)=\dim_\C(\Omega^{n+1}/\dd f\wedge\Omega^n)
\]
of $f$ is finite.
Then the following theorem holds.

\begin{thm}[C.~Sabbah {\cite[10.1--3]{Sab98b}}]\label{13}
$G_0$ is a free $\Ct$-- and $\Cth$--module of rank $\mu$.
\end{thm}

In particular, $G$ is a free $\Ctt$--module of rank $\mu$.

\section{Brieskorn lattice}

A finite $\Cta$-- resp. $\Cth$--submodule $L\subset G$ such that $L[\theta]=G$ resp. $L[\tau]=G$ is called a $\Cta$-- resp. $\Cth$--lattice.
By Theorem \ref{13}, a lattice is free of rank $\mu$.
In terms of a $\Cth$--basis of $G_0$, the $\Cthd$--module structure of $G$ is determined by the basis representation of $t$ on $G_0$.
The following lemma shows that the latter is determined by a matrix with coefficients in $\Cth$.

\begin{defn}
Let $\ul\phi$ be a basis of a $t$--invariant $\Cth$--sublattice $L\subset G$.
Then the matrix $A^\ul\phi\in\Cth^{\mu\times\mu}$ of $t$ with respect to $\ul\phi$ is defined by
\[
\ul\phi A^\ul\phi:=t\ul\phi.
\]
\end{defn}

\begin{lem}\label{21}
Let $\ul\phi$ be a basis of a $t$--invariant $\Cth$--sublattice $L\subset G$.
Then
\[
t\circ\ul\phi=\ul\phi\circ(A^\ul\phi+\theta^2\dth).
\]
\end{lem}
\begin{pf}
Since $t=\theta^2\dth$,
\begin{align*}
t\circ\ul\phi\Bigl(\sum_k\ol p_k\theta^k\Bigr)
&=t\sum_k\ul\phi\ol p_k\theta^k\\
&=\sum_kt(\ul\phi\ol p_k)\theta^k+\ul\phi\ol p_k\theta^2\dth\theta^k\\
&=\ul\phi\circ\bigl(A^\ul\phi+\theta^2\dth\bigr)\Bigl(\sum_k\ol p_k\theta^k\Bigr)
\end{align*}
and hence $t\circ\ul\phi=\ul\phi\circ(A^\ul\phi+\theta^2\dth)$.
\end{pf}

The following lemma gives a presentation of the $\Ctht$--module $G_0$.
This presentation shall be used to compute $t$ on $G_0$.

\begin{lem}\label{14}
There is an isomorphism of $\Ctht$--modules
\[
G_0=\Omega^{n+1}[\theta]/(\dd f-\theta\dd)(\Omega^n[\theta]).
\]
\end{lem}
\begin{pf}
By Theorem \ref{13},
\[
G=\Omega^{n+1}[\tau,\theta]/(\dd f-\theta\dd)(\Omega^n[\tau,\theta]).
\]
By definition, $G_0$ is the image of $\Omega^{n+1}[\theta]$ in $G$ and hence
\[
G_0=\Omega^{n+1}[\theta]\big/\bigl((\dd f-\theta\dd)(\Omega^n[\tau,\theta])\cap\Omega^{n+1}[\theta]\bigr).
\]
By Lemma \ref{39}, $\dd\ker(\dd f)\subset\dd f\wedge\dd\Omega^{n-1}\subset\ker(\dd f)$ and hence
\[
(\dd f-\theta\dd)(\Omega^n[\tau,\theta])\cap\Omega^{n+1}[\theta]=(\dd f-\theta\dd)(\Omega^n[\theta]).
\]
\end{pf}

Let $\ul x=x_0,\dots,x_n$ be coordinates on $\C^{n+1}$ with corresponding partial derivatives $\upd=\pd[x_0],\dots,\pd[x_n]$.
Let
\[
t:=f+\theta^2\dth\in\CXth\langle\dth\rangle.
\]
Then, by Lemma \ref{14}, we can identify
\[
G=\CXtt/(\upd(f)-\theta\upd)(\CXtt^{n+1})
\]
as $\Cttt$--modules and
\[
G_0=\CXth/(\upd(f)-\theta\upd)(\CXth^{n+1})
\]
as $\Ctht$--modules.
These modules are quotients of non--finite $\Cth$--modules.
On the numerator and denominator, the degree with respect to $\ul x$ defines an increasing filtration by finite $\Cth$--modules.
The following algorithm computes $t$ on $G_0$ by an approximation process with respect to these filtrations.

\begin{defn}
The degree $\deg_\ul x$ with respect to $\ul x$ defines an increasing filtration $\CXth_\bullet$ on $\CXth$ by finite $\Cth$--modules
\[
\CXth_k:=\{p\in\CXth\mid\deg_\ul x(p)\le k\}
\]
such that $t\CXth_\bullet\subset\CXth_{\bullet+\deg(f)}$.
We define the finite $\Cth$--modules
\begin{align*}
G_0^k&:=\CXth_k\big/\bigl((\upd(f)-\upd\theta)(\CXth^{n+1})\cap\CXth_k\bigr),\\
G_0^{k,l}&:=\CXth_k\big/\bigl((\upd(f)-\upd\theta)(\CXth_l^{n+1})\cap\CXth_k\bigr).
\end{align*}
\end{defn}

\begin{alg}\label{8}\
\begin{enumerate}
\item[Input:]
\begin{enumerate}
\item A cohomologically tame polynomial $f\in\CX$.
\item An integer $k\ge0$.
\end{enumerate}
\item[Output:]
\begin{enumerate}
\item A vector $\ul\phi\in\CXth^\mu$ such that $[\ul\phi]$ is a basis of a $t$--invariant $\Cth$--lattice $L_k\subset G_0$ and $L_k=G_0$ for $k\gg0$.
\item The matrix $A=A^{[\ul\phi]}\in\Cth^{\mu\times\mu}$.
\end{enumerate}
\item Set $l:=k$.
\item\label{11} Set $l:=l+1$.
\item Compute a reduced Gr\"obner basis
\[
\ul g:=\GB\bigl((\upd(f)-\upd\theta)(\ul x^\ul\alpha\ol e_i)\mid(\ul\alpha,i)\in\N^{n+1}\times[0,n],|\ul\alpha|\le l\bigr)
\]
of $(\upd(f)-\upd\theta)\bigl(\CXth_l^{n+1}\bigr)$ with respect to a monomial ordering $>$ on $\{\ul x^\ul\alpha\theta^i\mid(\ul\alpha,i)\in\N^{n+1}\times\N\}$ such that
\begin{align*}
|\ul\alpha|>|\ul\beta|&\Rightarrow\ul\alpha>\ul\beta,\\
(\ul\alpha,i)>(\ul\beta,j)&\Leftrightarrow\ul\alpha>\ul\beta\lor(\ul\alpha=\ul\beta\land i>j)
\end{align*}
for all $(\ul\alpha,i),(\ul\beta,j)\in\N^{n+1}\times\N$.
\item Find the minimal $k_0$ with
\[
k_0<|\alpha|\le k\Rightarrow\ul x^\ul\alpha\in\langle\lead(\ul g)\rangle\Cth.
\]
\item Compute $\ul\phi\in\CXth_{k_0}^\gamma$ such that $[\ul\phi]$ are cyclic generators of 
\[
G_0^{k_0,l}=\CXth_{k_0}/\langle g_i\mid\deg_\ul x(g_i)\le k_0\rangle\Cth\cong G_0^{k,l}
\]
and $\rho:=\rk(G_0^{k_0,l})$ using \cite[2.6.3]{GP02}.
\item\label{15} If $\rho>\mu$ or $\gamma>\rho=\mu$ then go to (\ref{11}).
\item\label{16} If $\rho<\mu$ then set $k:=k+1$ and go to (\ref{11}).
\item\label{17} If $k_0+\deg(f)>k$ then set $k:=k+1$ and go to (\ref{11}).
\item\label{18} If $[\ul\phi]$ is not a $\C$--basis of $\CX/\langle\upd(f)\rangle\CX$ then set $k:=k+1$ and go to (\ref{11}).
\item Compute a normal form $\NF(t\ul\phi,\ul g)$ of $t\ul\phi$ with respect to $\ul g$.
\item Compute the basis representation $A\in\Cth^{\mu\times\mu}$ of $[\NF(t\ul\phi,\ul g)]$ with respect to the $\Cth$--basis $[\ul\phi]$ of $L_k:=G_0^{k_0,l}$.
\item Return $\ul\phi$ and $A$.
\end{enumerate}
\end{alg}

\begin{lem}\label{28}
Algorithm \ref{8} terminates and is correct.
\end{lem}
\begin{pf}
Since $\upd(f)-\upd\theta$ is $\Cth$--linear,
\begin{multline*}
(\upd(f)-\upd\theta)(\CXth_l^{n+1})=\\\langle(\upd(f)-\upd\theta)(\ul x^\ul\alpha\ol e_i)\mid(\ul\alpha,i)\in\N^{n+1}\times[0,n],|\ul\alpha|\le l\rangle\Cth.
\end{multline*}
By definition of the monomial ordering,
\[
(\upd(f)-\upd\theta)(\CXth_l^{n+1})\cap\CXth_k=\langle g_i\mid\deg_\ul x(g_i)\le k\rangle\Cth
\]
and hence, by definition of $k_0$,
\begin{align*}
G_0^{k,l}&=\CXth_k/\langle g_i\mid\deg_\ul x(g_i)\le k\rangle\Cth\\
&\cong\CXth_{k_0}/\langle g_i\mid\deg_\ul x(g_i)\le k_0\rangle\Cth=G_0^{k_0,l}.
\end{align*}
Because of step (\ref{11}), $l$ is strictly increasing for fixed $k$.
There are $\Cth$--linear maps
\[
\xymat{
G_0^{k,l}\ar@{->>}[r]^-{\pi_{k,l}}&G_0^k\ar@{^(->}[r]^-{\iota_k}&G_0
}
\]
where $\iota_k$ is an isomorphism for $k\gg0$ and $\pi_{k,l}$ is an isomorphism for fixed $k$ and $l\gg0$.
By Theorem \ref{13}, $G_0^k$ is a free $\Cth$--module of rank at most $\mu$.
Hence, if condition (\ref{15}) holds then $\pi^{k,l}$ is not an isomorphism and if condition (\ref{16}) holds then $\iota^k$ is not an isomorphism.

By Theorem \ref{13}, there is a $\ul\psi\in\CXth^\mu$ such that $[\ul\psi]$ is a $\Cth$--basis of $G_0$.
In particular, $\iota^k$ is an isomorphism for $k\ge\deg_\ul x(\ul\psi)$ and hence
\[
\xymat{
\iota^k\circ\pi^{k,l}:G_0^{k,l}\ar@{^(->>}[r]&G_0
}
\]
is an isomorphism and conditions (\ref{15}) and (\ref{16}) do not hold for $k\ge\deg_\ul x(\ul\psi)$ and $l\gg0$.
By Lemma \ref{14}, for each $\ul\alpha\in\N^{n+1}$, there is a matrix $M^\ul\alpha\in\Cth^{\mu\times\mu}$ such that
\[
\ul x^\ul\alpha-\ul\psi M^\ul\alpha\in(\upd(f)-\upd\theta)(\CXth^{n+1}).
\]
If $|\ul\alpha|>\deg_\ul x(\ul\psi)$ then $\ul x^\ul\alpha-\ul\psi M^\ul\alpha\in(\upd(f)-\upd\theta)(\CXth_l^n)\cap\CXth_{|\ul\alpha|}$ and hence
\[
\ul x^\ul\alpha\in\lead((\upd(f)-\upd\theta)(\CXth_l^{n+1}))=\langle\lead(\ul g)\rangle\Cth
\]
for $l\gg0$.
Hence, by definition of $k_0$, $k_0\le\deg_\ul x(\ul\psi)$ for $k>\deg_\ul x(\ul\psi)$ and $l\gg0$ and, in particular, condition (\ref{17}) does not hold for $k\ge\deg_\ul x(\ul\psi)+\deg(f)$ and $l\gg0$.
Since $[\ul\psi]$ is a $\Cth$--basis of $G_0$, $[\ul\psi]$ is a $\C$--basis of
\[
G_0/\theta G_0=\CX/\langle\upd(f)\rangle\CX
\]
and hence condition (\ref{18}) does not hold for $k\ge\deg_\ul x(\ul\psi)$ and $l\gg0$.
This proves that the algorithm terminates.

Since $[\ul\phi]$ is a $\C$--basis of $\CX/\langle\upd(f)\rangle\CX=G_0/\theta G_0$, $\iota^k\circ\pi^{k,l}$ is injective and $[\ul\phi]$ is a basis of the $\Cth$--lattice $L_k=G_0^{k_0,l}\subset G_0$.
Since $\ul\phi\in\CXth_{k_0}^\mu$ and $k_0+\deg(f)\le k$, $t\ul\phi\in\CXth_k^\mu$ and hence, by definition of $k_0$, $\NF(t\ul\phi,\ul g)\in\CXth_{k_0}^\mu$.
By Lemma \ref{14},
\[
t[\ul\phi]=[t\ul\phi]=[\NF(t\ul\phi,\ul g)]=[\ul\phi A]=[\ul\phi]A
\]
and hence $L_k$ is $t$--invariant and $A=A^{[\ul\phi]}$.
This proves that the algorithm is correct.
\end{pf}

A priori, we do not know a $k_0$ such that $L_k=G_0$ for all $k\ge k_0$.
We shall solve this problem by a criterion on the spectrum with respect to the V--filtration.

\section{V--filtration}

\begin{defn}
The V--filtration $V_\bullet$ on $\Ctad$ is the increasing filtration by $V_0\Ctad$--modules
\begin{align*}
V_{-k}\Ctad&:=\tau^k\Cta\langle\tau\dta\rangle,\\
V_{k+1}\Ctad&:=V_k\Ctad+\dta V_k\Ctad
\end{align*}
for all $k\ge0$.
\end{defn}

\begin{prop}\label{20}
There is a unique $V_\bullet\Ctad$--good filtration $V_\bullet$ on $G$ by $\Cta$--lattices such that $\tau\dta+\alpha$ is nilpotent on $\gr^V_\alpha G$ for all $\alpha$.
\end{prop}
\begin{pf}
Since $G$ is regular at $0$, this follows from \cite[2.3.2, 4.1, 5.1.5]{Sab87}.
\end{pf}

\begin{defn}
$V_\bullet G$ is called the V--filtration on $G$.
\end{defn}

The following criterion shall be used to compute the V--filtration on $G$.

\begin{lem}\label{3}
Let $L\subset G$ be a $\tau\dta$--invariant $\Cta$--lattice with
\[
\spec(-\tau\dta\in\End(L/\tau L))\subset[\alpha,\alpha-1)
\]
for some $\alpha$. Then $L=V_\alpha G$.
\end{lem}
\begin{pf}
Let $\spec(-\tau\dta\in\End(L/\tau L))=\{\ul\alpha\}$ with 
\[
\alpha\ge\alpha_1>\cdots>\alpha_\nu>\alpha-1
\]
Let $\xymat{\phi:L/\tau L\ar[r]&L}$ be a $\Cta$--basis of $L$ and
\[
C_{\alpha_i}:=\phi(\ker((\tau\dta+\alpha_i)^\mu\in\End(L/\tau L)))
\]
for all $i\in[1,\nu]$.
Let 
\[
U_{\alpha_j-p}:=\tau^p\bigoplus_{i=j}^\nu C_{\alpha_i}\oplus\tau^{p+1}L
\]
for all $i\in[1,\nu]$ and $p\in\Z$.
Then $U_\bullet$ is an increasing filtration on $G$ by $\tau\dta$--invariant $\Cta$--lattices. 
By construction, $\tau\dta+\alpha_i-p$ is nilpotent on $\gr^U_{\alpha_i-p}G$ and $U_{\alpha_i-p}=\tau^pU_{\alpha_i}$ for all $i\in[1,\nu]$ and $p\in\Z$.
Since
\begin{align*}
\dta U_{\alpha_j-p}&=\tau^{p-1}(\tau\dta+p-1)\bigoplus_{i=j}^\nu C_{\alpha_i}\oplus\tau^p(\tau\dta+p)L\\
&\subset\tau^{p-1}(\tau\dta+p-1)\bigoplus_{i=j}^\nu C_{\alpha_i}\oplus\tau^p(\tau\dta+p)\bigoplus_{i=1}^{j-1}C_{\alpha_i}+U_{\alpha_j-p},\\
U_{\alpha_j-p+1}&=\tau^{p-1}\bigoplus_{i=j}^\nu C_{\alpha_i}\oplus\tau^p\bigoplus_{i=1}^{j-1}C_{\alpha_i}\oplus\tau^{p+1}L\\
&\subset\tau^{p-1}\bigoplus_{i=j}^\nu C_{\alpha_i}\oplus\tau^p\bigoplus_{i=1}^{j-1}C_{\alpha_i}+U_{\alpha_j-p},
\end{align*}
$\dta U_{\alpha_j-p}+U_{\alpha_j-p}=U_{\alpha_j-p+1}$ for $p>\alpha_j+1$ and hence $U_\bullet$ is $V_\bullet\Ctad$--good.
Then, by Proposition \ref{20}, $U_\bullet G=V_\bullet G$ and hence $L=V_\alpha G$.
\end{pf}

The following algorithm computes the V--filtration using the criterion in Lemma \ref{3}.
For a given $\Cth$--lattice with $\Cth$--basis $\ul\phi$, $L:=\langle\ul\phi\rangle\Cta$ is a $\Cta$--lattice with $\Cta$-basis $\ul\phi$ and $-\tau\dta\ul\phi=\ul\phi B$ where $B=\tau A^\ul\phi\in\Ctt^{\mu\times\mu}$.
By a saturation process of $L$ with respect to $\tau\dta$, $L$ is replaced by a $\tau\dta$--invariant $\Cta$--lattice and $\ul\phi$ is modified such that $B\in\Cta^{\mu\times\mu}$.
Then a sequence of basis transformations modifies $\ul\phi$ such that $\spec(B_0)\subset[\alpha,\alpha-1)$ for some $\alpha$.

\begin{alg}\label{1}\
\begin{enumerate}
\item[Input:] The matrix $A=A^\ul\phi\in\Cth^{\mu\times\mu}$ for a basis $\ul\phi$ of a $t$--invariant $\Cth$--lattice $L\subset G$.
\item[Output:]
\begin{enumerate}
\item A matrix $U\in\Cth^{\mu\times\mu}$ such that $\ul\phi U$ is a $\Cta$--basis of $V_\alpha$ for some $\alpha$.
\item A matrix $B=\sum_{i\ge0}B_i\tau^i\in\Cta^{\mu\times\mu}$ such that $-\tau\dta(\ul\phi U)=\ul\phi U B$ and $\spec(B_0)=\{\ul\alpha\}$ with $\alpha\ge\alpha_1>\cdots>\alpha_\nu>\alpha-1$.
\end{enumerate}
\item
\begin{enumerate}
\item Set $k:=0$ and $U_0:=E\in\C^{\mu\times\mu}$.
\item Until $\{(\tau A-\tau\dta)(U_k)\}\subset\langle U_k\rangle\Cta$ do:
\begin{enumerate}
\item Set $k:=k+1$.
\item Compute $U_k\in\Cth^{\mu\times\mu}$ with $\deg(U_k)\le k(\deg(A)-1)$ such that
\[
\langle U_{k+1}\rangle\Cta=\langle U_k\rangle\Cta+\langle(\tau A-\tau\dta)(U_k)\rangle\Cta.
\]
\end{enumerate}
\item Set $U:=U_k$.
\end{enumerate}
\item
\begin{enumerate}
\item Set $B=\sum_{i\ge0}B_i\tau^i:=\inv U(\tau A-\tau\dta)(U)\in\Cta^{\mu\times\mu}$.
\item Compute $\{\ul\alpha\}:=\spec(B_0)$ and $j\in[1,\nu]$ such that
\[
\alpha_1>\dots>\alpha_j>\alpha_1-1\ge\alpha_{j+1}>\dots>\alpha_\nu.
\]
\item\label{2} If $j=\nu$ then return $U$ and $B$.
\item Compute $U_0\in\GL_\mu(\C)$ such that
\[
\inv U_0BU_0=
\begin{pmatrix}
B^{1,1}&B^{1,2}\\
B^{2,1}&B^{2,2}
\end{pmatrix}
\]
where
\begin{align*}
\spec(B^{1,1}_0)&=\{\alpha_1,\dots,\alpha_j\},\\
\spec(B^{2,2}_0)&=\{\alpha_{j+1},\dots,\alpha_\nu\},
\end{align*}
$B^{1,2}_0=0$, and $B^{2,1}_0=0$.
\item Set $U=\begin{pmatrix}U_1&U_2\end{pmatrix}:=UU_0$ and
\[
B=
\begin{pmatrix}
B^{1,1}&B^{1,2}\\
B^{2,1}&B^{2,2}
\end{pmatrix}
:=\inv U_0BU_0.
\]
\item Set $U:=\begin{pmatrix}U_1&\inv\tau U_2\end{pmatrix}$ and
\[
B:=\begin{pmatrix}
B^{1,1}&\inv\tau B^{1,2}\\
\tau B^{2,1}&B^{2,2}+E
\end{pmatrix}.
\]
\item Set $\alpha_i:=\alpha_i+1$ for $i=j+1,\dots,\nu$.
\item Reorder $\ul\alpha$ and redefine $j\in[1,\nu]$ such that
\[
\alpha_1>\dots>\alpha_j>\alpha_1-1\ge\alpha_{j+1}>\dots>\alpha_\nu.
\]
\item Go to (\ref{2}).
\end{enumerate}
\end{enumerate}
\end{alg}

\begin{rem}\
\begin{enumerate}
\item If $A=A_0+\theta A_1$ then $U_k=U_0=E$.
\item If
\[
B=\begin{pmatrix}
B^{1,1}&B^{1,2}\\
0&B^{2,2}
\end{pmatrix}
\]
with $\spec(B^{1,1}_0)=\{\alpha_1,\dots,\alpha_j\}$ and $\spec(B^{2,2}_0)=\{\alpha_{j+1},\dots,\alpha_\nu\}$ then one can choose
\[
U_0=\begin{pmatrix}
E&U_0^{1,2}\\
0&E
\end{pmatrix}.
\]
\end{enumerate}
\end{rem}

\begin{lem}\label{29}
Algorithm \ref{1} terminates and is correct.
\end{lem}
\begin{pf}\
\begin{enumerate}
\item By Lemma \ref{21},
\begin{align*}
\langle\ul\phi U_{k+1}\rangle\Cta&=\langle\ul\phi U_k\rangle\Cta+\langle\ul\phi\circ(\tau A-\tau\dta)(U_k)\rangle\Cta\\&=\langle\ul\phi U_k\rangle\Cta+\tau\langle\ul\phi\circ(A+\theta^2\dth)(U_k)\rangle\Cta\\&=\langle\ul\phi U_k\rangle\Cta+\tau\dta\langle\ul\phi U_k\rangle\Cta
\end{align*}
and hence $\{\langle\ul\phi U_k\rangle\Cta\}_{k\ge0}$ is an increasing sequence of finite $\Cta$--modules.
Since $\langle\ul\phi U_0\rangle\Cta=\Cta^\mu$, one can choose $U_k\in\Cth^{\mu\times\mu}$.
The V--filtration on $G$ consists of finite and hence Noetherian $\tau\dta$--invariant $\Cta$--modules.
For some $\alpha$, $\{\ul\phi\}\subset V_\alpha G$ and hence $\langle\ul\phi U_k\rangle\Cta\subset V_\alpha$ for all $k\ge0$.
This implies that the sequence $\{\langle U_k\rangle\Cta\}_{k\ge0}$ is stationary.
Then $\langle\ul\phi U\rangle\Cta\subset G$ is a $\tau\dta$--invariant $\Cta$--lattice.
\item By Lemma \ref{21}, $-\tau\dta\circ\ul\phi=\tau t\circ\ul\phi=\ul\phi\circ(\tau A-\tau\dta)$ and hence
\[
-\tau\dta\circ\ul\phi U=\ul\phi U\circ(B-\tau\dta).
\]
The $\tau\dta$--invariance of the $\Cta$--lattice $\langle\ul\phi U\rangle\Cta$ is preserved since
\begin{align*}
\begin{pmatrix}U_1&\inv\tau U_2\end{pmatrix}
&=
\begin{pmatrix}U_1&U_2\end{pmatrix}
\begin{pmatrix}
E&\\
&\inv\tau E
\end{pmatrix},\\
\begin{pmatrix}
B^{1,1}&\inv\tau B^{1,2},\\
\tau B^{2,1}&B^{2,2}+E
\end{pmatrix}
&=
\begin{pmatrix}
E&\\
&\tau E
\end{pmatrix}
\biggl(
\begin{pmatrix}
B^{1,1}&B^{1,2}\\
B^{2,1}&B^{2,2}
\end{pmatrix}
-\tau\dta
\biggr)
\begin{pmatrix}
E&\\
&\inv\tau E
\end{pmatrix},
\end{align*}
and $B^{1,2}\inv\tau\in\Cta^{j\times(\mu-j)}$.
The index $j$ is strictly increasing since
\begin{align*}
\spec(B_0)&=\spec\begin{pmatrix}
B^{1,1}_0-E&\inv\tau B^{1,2}_0\\
0&B^{2,2}_0
\end{pmatrix}\\
&=\{\alpha_1,\dots,\alpha_j,\alpha_{j+1}+1,\dots,\alpha_\nu+1\}
\end{align*}
and hence the algorithm terminates.
Then $L:=\langle\ul\phi U\rangle\Cta\subset G$ is a $\tau\dta$--invariant $\Cta$--lattice with $\spec(-\tau\dta\in\End(L/\tau L))\subset[\alpha,\alpha-1)$ for $\alpha:=\alpha_1$.
Hence, by Lemma \ref{3}, $L=V_\alpha$ and $\ul\phi U$ is a $\Cta$--basis of $V_\alpha G$.
\end{enumerate}
\end{pf}

\section{Spectrum}

The spectrum with respect to the V--filtration shall be used to check equality of $\Cth$--lattices.

\begin{defn}\
\begin{enumerate}
\item The spectrum $\xymat{\spec(F_\bullet):\Q\ar[r]&\N}$ of an increasing filtration $F_\bullet$ on a finite vector space $V$ is defined by
\[
\spec(F_\bullet)(\alpha):=\dim(\gr^{F_\bullet}_\alpha V)
\]
for all $\alpha\in\Q$.
The spectrum $\spec(F_\bullet)$ of a decreasing filtration $F^\bullet$ on $V$ is defined analogously.
\item The spectrum of a $\Cth$--lattice $L\subset G$ is defined by
\[
\spec(L):=\spec(V_\bullet(L/\theta L)).
\]
\item The spectrum of $f$ is defined by
\[
\spec(f):=\spec(G_0).
\]
\end{enumerate}
\end{defn}

The following algorithm computes the spectrum of a $t$--invariant $\Cth$--lattice by computing a Gr\"obner basis compatible with the V--filtration.

\begin{alg}\label{22}\
\begin{enumerate}
\item[Input:]
\begin{enumerate}
\item A matrix $B=\sum_{i\ge0}B_i\tau^i\in\Cta^{\mu\times\mu}$ such that $-\tau\dta\ul\phi=\ul\phi B$ for a $\Cta$--basis $\ul\phi$ of $V_\alpha$ and $\spec(B_0)=\{\ul\alpha\}$ with $\alpha\ge\alpha_1>\cdots>\alpha_\nu>\alpha-1$.
\item A matrix $M\in\Ctt^{\mu\times\mu}$ such that $\ul\phi M$ is a basis of a $t$--invariant $\Cth$--lattice $L\subset G$.
\end{enumerate}
\item[Output:] The spectrum $\sigma=\spec(L)\in\Q^\N$.
\item Compute $U_0\in\GL_\mu(\C)$ such that
\[
\inv U_0B_0U_0=
\begin{pmatrix}
B_0^1&&\\
&\ddots&\\
&&B_0^\nu
\end{pmatrix}
\]
where $B_0^i\in\C^{\mu_i\times\mu_i}$ with $\spec(B_0^i)=\{\alpha_i\}$ for $i\in[1,\nu]$.
\item Set $\ul\phi=(\ul\phi^i)_{i\in[1,\nu]}:=\ul\phi U_0$ and $M:=\inv U_0M$.
\item Compute a minimal Gr\"obner basis
\[
M:=\GB(M)\in\Cth^{\mu\times\mu}
\]
compatible with the ordering $>$ on $\{\theta^k\ul\phi^i\mid(k,i)\in\Z\times[1,\nu]\}$ defined by
\[
(k,i)>(l,j):\Leftrightarrow k>l\lor(k=l\land i>j)
\]
for all $(k,i),(l,j)\in\Z\times[1,\nu]$.
\item Return $\sigma\in\Q^\N$ with
\[
\sigma(k+\alpha_i):=\#\bigl(\{\lead(M)\}_{(k,i)}\bigr)
\]
for all $(k,i)\in\Z\times[1,\nu]$.
\end{enumerate}
\end{alg}

\begin{lem}\label{30}
Algorithm \ref{22} terminates and is correct.
\end{lem}
\begin{pf}
Since $-\tau\dta\ul\phi=\ul\phi B$,
\[
-\tau\dta(\theta^l\ul\phi^j)\equiv\theta^q\ul\phi^j(B_0^j+q)\mod\bigoplus_{(k,i)\le(l,j)}\langle\theta^k\ul\phi^i\rangle\C
\]
with $\spec(B_0^j+l)=\{\alpha_j+l\}$.
Then, by Lemma \ref{3},
\[
V_{\alpha_j+l}G=\bigoplus_{(i,k)\le(j,l)}\theta^k\langle\ul\phi^i\rangle\C
\]
and hence, since $M$ is a minimal Gr\"obner basis, 
\begin{align*}
\spec(L)(\alpha_j+l)&=\dim_\C\gr^V_{\alpha+l}(L/\theta L)\\
&=\dim_\C\bigl((\gr^VL/\theta\gr^VL)_{\alpha_j+l}\bigl)\\
&=\dim_\C\bigl(\langle\lead(M)\rangle\Cth/\theta\langle\lead(M)\rangle\Cth)_{(q,j)}\bigr)\\
&=\dim_\C\bigl((\langle\lead(M)\rangle\C)_{(l,j)}\bigr)\\
&=\#\bigl(\{\lead(M)\}_{(l,j)}\bigr)
\end{align*}
for all $(l,j)\in\Z\times[1,\nu]$.
\end{pf}

The following lemma reduces the problem of equality of $\Cth$--lattices to the problem of equality of filtrations on a finite vector space.

\begin{defn}\label{35}
A $\Cth$--lattice $L\subset G$ defines an increasing filtration $L_\bullet$ on $G$ by $\Cth$--lattices $L_p:=\tau^pL$ and a corresponding decreasing filtration $L^\bullet:=L_{n-\bullet}$.
We denote the filtrations defined by $G_0$ by $G_\bullet$ and $G^\bullet$.
\end{defn}

\begin{lem}\label{6}
Let $L\subset G$ be a $\Cth$--lattice.
Then
\[
\xymat{\gr_L^{n-p}\gr^V_\alpha G\ar[r]<2pt>^-{\theta^p}&\ar[l]<2pt>^-{\tau^p}\gr^V_{\alpha+p}\gr_L^nG=\gr^V_{\alpha+p}(L/\theta L)}
\]
is an isomorphism for all $\alpha\in\Q$ and $p\in\Z$.
\end{lem}
\begin{pf}
This follows from $\theta^pV_\alpha G=V_{\alpha+p}G$ and $L^{n-p}=\tau^pL$ for all $\alpha\in\Q$ and $p\in\Z$.
\end{pf}

\begin{defn}
The sum $\xymat{\sum:\N^\Q\ar[r]&\Q}$ is defined by
\[
\sum\sigma:=\sum_{\alpha\in\Q}\alpha\sigma(\alpha)
\]
for all $\sigma\in\N^\Q$.
\end{defn}

The following lemma gives a criterion on the mean value of the spectrum to check equality of filtrations on a finite vector space.

\begin{lem}\label{10}
Let $F_1^\bullet$ and $F_2^\bullet$ be decreasing filtrations on a finite vector space $V$ with $F_2^\bullet\subset F_1^\bullet$.
Then $\sum\spec(F_2^\bullet)\le\sum\spec(F_1^\bullet)$ and equality implies that $F_2^\bullet=F_1^\bullet$.
\end{lem}
\begin{pf}
This is an elementary fact from linear algebra.
\end{pf}

The following criterion on the mean value of the spectrum shall be used to check equality of $\Cth$--lattices.

\begin{lem}\label{27}
Let $L_2\subset L_1\subset G$ be $\Cth$--lattices. 
Then $\sum\spec(L_1)\le\sum\spec(L_2)$ and equality implies that $L_1=L_2$.
\end{lem}
\begin{pf}
Since $L_2\subset L_1$,
\[
L_2^\bullet\gr^V_{[0,1)}G\subset L_1^\bullet\gr^V_{[0,1)}G
\]
where $\gr^V_{[0,1)}G=\bigoplus_{0\le\alpha<1}\gr^V_\alpha G$ and hence, by Lemma \ref{10},
\[
\sum\spec\bigl(L_2^\bullet\gr^V_{[0,1)}G\bigr)\le\sum\spec\bigl(L_1^\bullet\gr^V_{[0,1)}G\bigr).
\]
By Lemma \ref{6}, $\spec(L_i)(\alpha+p)=\spec\bigl(L_i^\bullet\gr^V_\alpha G\bigr)(n-p)$ and hence
\begin{align*}
\sum\spec(L_i)&=\sum_{0\le\alpha<1}\sum_{p\in\Z}(\alpha+n-p)\spec\bigl(L_i^\bullet\gr^V_\alpha G\bigr)(p)\\
&=n\mu+\sum_{0\le\alpha<1}\alpha\dim_\C\bigl(\gr^V_\alpha G\bigr)-\sum\spec\bigl(L_i^\bullet\gr^V_\alpha G\bigr)\\
&=n\mu+\sum_{0\le\alpha<1}\alpha\dim_\C\bigl(\gr^V_\alpha G\bigr)-\sum\spec\bigl(L_i^\bullet\gr^V_{[0,1)}G\bigr)
\end{align*}
for $i=1,2$.
This implies that
\[
\sum\spec(L_2)-\sum\spec(L_1)=\sum\spec\bigl(L_1^\bullet\gr^V_{[0,1)}G\bigr)-\sum\spec\bigl(L_2^\bullet\gr^V_{[0,1)}G\bigr).
\]
Let $x\in(L_1\backslash L_2)\cap\bigl(V_{\alpha+p}G\big\backslash V_{<\alpha+p}G\bigr)$ with $0\le\alpha<1$ and minimal $\alpha+p$.
Then, in particular, $x\notin\theta L_1$ and hence, by Lemma \ref{6}, $0\neq[\tau^p x]\in\gr_{L_1}^{n-p}\gr^V_\alpha G$.
Moreover, there is a $q\ge1$ such that $\theta^qx\in L_2\backslash\theta L_2$ and, again by Lemma \ref{6}, $0\ne[\tau^p x]\in\gr_{L_2}^{n-p-q}\gr^V_\alpha G$.
This implies that $L_2^{n-p}\gr^V_\alpha G\subsetneq L_1^{n-p}\gr^V_\alpha G$ and hence
\[
L_2^\bullet\gr^V_{[0,1)}G\subsetneq L_1^\bullet\gr^V_{[0,1)}G.
\]
Then the claim follows from Lemma \ref{10}.
\end{pf}

The following theorem gives the mean value of the spectrum of $G_0$.

\begin{thm}[C.~Sabbah {\cite[11.1]{Sab98b}}]\label{26}
$\frac{1}{\mu}\sum\spec(G_0)=\frac{n+1}{2}$.
\end{thm}

By Theorem \ref{26}, one can compute $t$ on $G_0$ using Algorithm \ref{8}, \ref{1}, and \ref{22} by increasing $k$ until $\frac{1}{\mu}\sum\spec(L_k)=\frac{n+1}{2}$.

Our final goal is to compute a good basis of $G_0$.
In terms of a good basis of $G_0$, the matrix of $t$ has degree one and its degree one part determines the spectrum of $f$.

\begin{defn}
Let $\ul\phi$ be a $\Cth$--basis of a $t$--invariant $\Cth$--lattice $L\subset G$.
Then $\ul\phi$ is called good if $A^{\ul\phi}=A^{\ul\phi}_0+\theta A^{\ul\phi}_1$ where $A^\ul\phi_0,A^\ul\phi_1\in\C^{\mu\times\mu}$,
\[
A^{\ul\phi}_1=\begin{pmatrix}
\alpha_1 & & \\
 & \ddots & \\
 & & \alpha_\mu
\end{pmatrix}
\]
and $\phi_i\in V_{\alpha_i}L$ for all $i\in[1,\mu]$.
\end{defn}

\begin{lem}\label{37}
Let $\ul\phi$ be a good basis of a $t$-invariant $\Cth$--lattice $L\subset G$ and
\[
\begin{pmatrix}
\alpha_1 & & \\
 & \ddots & \\
 & & \alpha_\mu
\end{pmatrix}:=A^\ul\phi_1.
\]
Then $\spec(L)(\alpha)=\#\{i\in[1,\mu]\mid\alpha_i=\alpha\}$.
\end{lem}
\begin{pf}
Since $\ul\phi$ is a $\Cth$--basis of $L$ and $\phi_i\in V_{\alpha_i}$ for all $i\in[1,\mu]$, $\bigl([\phi_i]\bigr)_{\alpha_i=\alpha}$ is a $\C$--basis of $\gr^V_\alpha(L/\theta L)$ and hence $\spec(L)(\alpha)=\#\{i\in[1,\mu]\mid\alpha_i=\alpha\}$.
\end{pf}

\section{Monodromy}

Let $T_\infty$ be the monodromy of $M$ around the discriminant $\dsc(f)=f(\crt(f))$ of $f$ and $\wh T_0$ be the monodromy of $G$ at $0$.

\begin{thm}[C.~Sabbah {\cite[1.10]{Sab98a}}]\label{33}
$T_\infty=\wh T_0^{-1}$.
\end{thm}

Using Theorem \ref{33}, the monodromy $T_\infty$ can be read off from the matrix of $t$ with respect to a good basis.

\begin{prop}\label{34}
Let $\ul\phi$ be a good basis of a $\Cth$--lattice $L\subset G$ and
\[
\begin{pmatrix}
\alpha_1 & & \\
 & \ddots & \\
 & & \alpha_\mu
\end{pmatrix}:=A^\ul\phi_1.
\]
Then
\[
\exp\bigl(-2\pi\ii\bigl(\gr^V_1\bigl(A^\ul\phi_0\bigr)+A^\ul\phi_1\bigr)\bigr)
\]
is a matrix of $T_\infty$ where
\[
(V_\alpha(C))_{i,j}:=\begin{cases}
c_{i,j} & \text{if } \alpha_i\ge\alpha_j+\alpha,\\
0 & \text{else,}
\end{cases}
\]
for $C=(c_{i,j})_{i,j}\in\C^{\mu\times\mu}$.
\end{prop}
\begin{pf}
Since $\ul\phi$ is a $\Cth$--basis of $L$ and $\phi_i\in V_{\alpha_1}$ for all $i\in[1,\mu]$, $\bigl([\phi_i]\bigr)_{\alpha_i=\alpha}$ is a $\C$--basis of $\gr^V_\alpha(L/\theta L)$ and hence, by Lemma \ref{6},
\[
\xymat{
-\tau\dta=\dt t:\gr_L^p\gr^V_\alpha G\ar[r]^-{\alpha\oplus N}&\gr_L^p\gr^V_\alpha G\oplus\gr_L^{p-1}\gr^V_\alpha G
}
\]
for all $\alpha\in\Q$ and $p\in\Z$ where $\theta N$ is induced by
\[
\gr^V_1\gr_L^0t=\gr^V_1\bigl(\ul\phi\circ A^\ul\phi_0\circ\ul\phi^{-1}\bigr)=\ul\phi\circ\gr^V_1\bigl(A^\ul\phi_0\bigr)\circ\ul\phi^{-1}.
\]
Then, by \cite[6.0.1]{Sab87}, $\exp\bigl(2\pi\ii\bigl(\gr^V_1\bigl(A^\ul\phi_0\bigr)+A^\ul\phi_1\bigr)\bigr)$ is a matrix of $\wh T_0$ and hence, by Theorem \ref{33}, $\exp\bigl(-2\pi\ii\bigl(\gr^V_1\bigl(A^\ul\phi_0\bigr)+A^\ul\phi_1\bigr)\bigr)$ is a matrix of $T_\infty$.
\end{pf}

\section{Good lattices}

The following property is sufficient for the existence of a good basis of a $\Cth$--lattice \cite[5.2]{Sab98b}.
Recall that a morphism $\xymat{N:F_1^\bullet V_1\ar[r]&F_2^\bullet V_2}$ of filtered vector spaces $F_i^\bullet V_i$ for $i=1,2$ is called strict if $N(V_1)\cap F_2^p=N(F_2^p)$ for all $p\in\Z$.

\begin{defn}
We call a $t$--invariant $\Cth$--lattice $L\subset G$ good if
\[
\xymat{(\tau\dta+\alpha)^p:L^\bullet\gr^V_\alpha G\ar[r]&L^{\bullet-p}\gr^V_\alpha G}
\]
is strict for all $\alpha\in\Q$ and $p\ge1$.
\end{defn}

The following theorem follows from the fact that
\[
N:=\bigoplus_{0\le\alpha<1}(\tau\dta+\alpha)
\]
is a morphism of a natural mixed Hodge structure on the moderate nearby cycles $\psi^\text{mod}_\tau G=\bigoplus_{0\le\alpha<1}\gr^V_\alpha G$ with Hodge filtration induced by $G_\bullet$ as defined in Definition \ref{35} \cite[13.1]{Sab98b}.

\begin{thm}[C.~Sabbah {\cite[13.3]{Sab98b}}]\label{23}
$G_0$ is a good $\Cth$--lattice.
\end{thm}

The following lemma shall be used to construct an opposite filtration of $L^\bullet$ on $\gr^VG$ for a good lattice $L$.

\begin{lem}\label{5}
Let $V$ be a finite vector space, $F^\bullet$ a decreasing filtration on $V$ with $F^p=0$ for $p>m$, and $N\in\End(V)$ such that
\[
\xymat{N^p:F^\bullet V\ar[r]&F^{\bullet-p}V}
\]
strict for all $p\ge1$.
Then $\sum_{q\ge0}N^q(F^m)=\bigoplus_{q\ge0}N^q(F^m)$ and
\[
\xymat{N^p:F^\bullet\bigl(V\big/\sum_{q\ge0}N^q(F^m)\bigr)\ar[r]&F^{\bullet-p}\bigl(V\big/\sum_{q\ge0}N^q(F^m)\bigr)}
\]
is strict for all $p\ge1$.
\end{lem}
\begin{pf}
If $x\in F^m$ with $N^{p+1}(x)\in\sum_{q=0}^pN^q(F^m)\subset F^{m-p}$ then
\[
N^{p+1}(x)\in N^{p+1}(F^{m+1})=0
\]
since $N^{p+1}$ is strict and $F^{m+1}=0$.
Hence,
\[
N^{p+1}(F^m)+\sum_{q=0}^pN^q(F^m)=N^{p+1}(F^m)\oplus\sum_{q=0}^pN^q(F^m)
\]
and, by induction, $\sum_{q\ge0}N^q(F^m)=\bigoplus_{q\ge0}N^q(F^m)$.

Let $N^p(x)\in F^q+\sum_{r\ge0}N^r(z_r)$ with $z_r\in F^m$ for all $r\ge0$.
If $m-p<q$ then $N^p\bigl(x-\sum_{r>p}N^{r-p}(z_r)\bigr)\in F^{m-p+1}$ and hence
\[
N^p(x)\in N^p(F^{m+1})+\sum_{r>p}N^r(F^m)\subset N^p(F^{q+p})+\sum_{r\ge0}N^r(F^m)
\]
since $N^p$ is strict and $F^{m+1}=0$.
If $m-p\ge q$ then $N^p\bigl(x-\sum_{r>p}N^{r-p}(z_r)\bigr)\in F^q$ and hence $N^p(x)\in N^p(F^{q+p})+\sum_{r\ge0}N^r(F^m)$ since $N^p$ is strict.
This implies that $N^p$ is strict modulo $\sum_{r\ge0}N^r(F^m)$ for all $p\ge1$.
\end{pf}

The following algorithm computes a $\Ctt$--basis of $G$ compatible with the V--filtration refined by an opposite filtration of $L^\bullet$ on $\gr^VG$ for a good lattice $L$.
This basis shall be used to compute a good basis of $L$.

\begin{alg}\label{4}\ 
\begin{enumerate}
\item[Input:]
\begin{enumerate}
\item A matrix $B=\sum_{i\ge0}B_i\tau^i\in\Cta^{\mu\times\mu}$ such that $-\tau\dta\ul\phi=\ul\phi B$ for a $\Cta$--basis $\ul\phi$ of $V_\alpha$ and $\spec(B_0)=\{\ul\alpha\}$ with $\alpha\ge\alpha_1>\cdots>\alpha_\nu>\alpha-1$.
\item A matrix $M\in\Ctt^{\mu\times\mu}$ such that $\ul\phi M$ is a basis of a good $\Cth$--lattice $L\subset G$.
\end{enumerate}
\item[Output:]
A matrix $U=(U^{i,p})_{(i,p)\in[1,\nu]\times\Z}\in\GL_\mu(\C)$ such that $(\ul\phi U^{i,q})_{q\ge p}$ is a $\C$--basis of $L^p\gr^V_{\alpha_i}G$ and $\{(\theta\dth-\alpha_i)(\ul\phi U^{i,p})\}\subset\{\ul\phi U^{i,p-1}\}+V_{\alpha-1}$ for all $(i,p)\in[1,\nu]\times\Z$.
\item Compute $U_0\in\GL_\mu(\C)$ such that
\[
\inv U_0B_0U_0=
\begin{pmatrix}
B_0^1&&\\
&\ddots&\\
&&B_0^\nu
\end{pmatrix}
\]
where $B_0^i\in\C^{\mu_i\times\mu_i}$ with $\spec(B_0^i)=\{\alpha_i\}$ for $i\in[1,\nu]$.
\item Set $\ul\phi=(\ul\phi^i)_{i\in[1,\nu]}:=\ul\phi U_0$ and $M:=\inv U_0M$.
\item Compute a Gr\"obner basis
\[
M:=\GB(M)\in\Cth^{\mu\times\mu}
\]
compatible with the ordering $>$ on $\{\theta^p\ul\phi^i\mid(p,i)\in\Z\times[1,\nu]\}$ defined by
\[
(p,i)>(q,j):\Leftrightarrow p>q\lor(p=q\land i>j)
\]
for all $(p,i),(q,j)\in\Z\times[1,\nu]$.
\item Set $(M^{p,i})_{(p,i)\in\Z\times[1,\nu]}:=M$ where
\[
\{\lexp(M^{p,i})\}=\{(p,i)\}
\]
for all $(p,i)\in\Z\times[1,\nu]$.
\item For $i=1,\dots,r$ do:
\begin{enumerate}
\item Compute $(F^{i,p})_{p\in\Z}\in\C^{\mu_i\times\mu_i}$ such that
\[
F^{i,p}:=(\tau^q\lead(M^{q,i}))_{q\le n-p}.
\]
\item Set $N_i:=B_0^i-\alpha_i$.
\item Compute $U^i=(U^{i,p})_{p\in\Z}\in\C^{\mu_i\times\mu_i}$ such that
\[
\langle F^{i,p}\rangle\C=\langle U^{i,q}\mid q\le p\rangle\C,\quad\{N_iU^{i,p}\}\subset\{U^{i,p-1}\}
\]
for all $p\in\Z$.
\end{enumerate}
\item Return
\[
U=(U^{i,p})_{(i,p)\in[1,\nu]\times\Z}:=U^0
\begin{pmatrix}
U^1&&\\
&\ddots&\\
&&U^\nu
\end{pmatrix}.
\]
\end{enumerate}
\end{alg}

\begin{lem}\label{31}
Algorithm \ref{4} terminates and is correct.
\end{lem}
\begin{pf}
Since $-\tau\dta\ul\phi=\ul\phi B$,
\[
-\tau\dta(\theta^q\ul\phi^j)\equiv\theta^q\ul\phi^j(B_0^j+q)\mod V_{q+\alpha-1}
\]
with $\spec(B_0^j+q)=\{\alpha_j+q\}$ and hence, by Lemma \ref{3},
\[
V_{\alpha_j+q}G=\bigoplus_{(p,i)\le(q,j)}\theta^p\langle\ul\phi^i\rangle\C
\]
for all $(q,j)\in\Z\times[1,\nu]$.
Since $M$ is a Gr\"obner basis, this implies that $\ul\phi M^{q,j}\in V_{\alpha_j+q}L$ for all $(q,j)\in\Z\times[1,\nu]$.
Then, by Lemma \ref{6},
\[
L^\bullet\gr^V_{\alpha_i}G=\bigl(\langle\ul\phi^iF^{i,\bullet}\rangle\C+V_{\alpha_i}\bigr)\big/V_{<\alpha_i}\subset\gr^V_{\alpha_i}G
\]
and, since $L$ is good and $(\theta\dth-\alpha_i)\ul\phi^i\equiv\ul\phi^iN_i\mod V_{\alpha-1}$,\[
\xymat{N_i^p:\langle F^{i,\bullet}\rangle\C\ar[r]&\langle F^{i,\bullet-p}\rangle\C}
\]
is strict for all $i\in[1,\nu]$ and $p\ge1$.
Hence, by Lemma \ref{5}, one can compute $U^i=(U^{i,p})_{p\in\Z}\in\C^{\mu_i\times\mu_i}$ such that
\[
\langle F^{i,p}\rangle\C=\langle U^{i,q}\mid q\ge p\rangle\C,\quad\{N_iU^{i,p}\}\subset\{U^{i,p-1}\}
\]
for all $(i,p)\in[1,\nu]\times\Z$.
Then
\[
(\theta\dth-\alpha_i)(\ul\phi U^{i,p})\equiv\ul\phi^iN_iU^{i,p}\mod V_{\alpha-1}
\]
and hence $\{(\theta\dth-\alpha_i)(\ul\phi U^{i,p})\}\subset\{\ul\phi U^{i,p-1}\}+V_{\alpha-1}$ for all $(i,p)\in[1,\nu]\times\Z$.
\end{pf}

\section{Good bases}

The following algorithm computes a good basis of a good lattice $L$ by a simultaneous normal form computation and basis transformation.
The computation requires a $\Ctt$--basis of $G$ compatible with the V--filtration refined by an opposite filtration of $L^\bullet$ on $\gr^VG$.

\begin{alg}\label{9}\
\begin{enumerate}
\item[Input:]
\begin{enumerate}
\item A matrix $B\in\Cta^{\mu\times\mu}$ with $\spec(B_0)=\{\ul\alpha\}$ and $\alpha\ge\alpha_1>\cdots>\alpha_\nu>\alpha-1$ such that $-\tau\dta\ul\phi=\ul\phi B$ for a $\Cta$--basis $\ul\phi$ of $V_\alpha$
\item A matrix $M\in\Ctt^{\mu\times\mu}$ such that $\ul\phi M$ is a basis of a good $\Cth$--lattice $L\subset G$.
\item An indexing $\ul\phi=(\ul\phi^{i,p})_{(i,p)\in[1,\nu]\times\Z}$ such that $(\ul\phi^{i,q})_{q\ge p}$ is a $\C$--basis of $L^p\gr^V_{\alpha_i}G$ and $\{(\theta\dth-\alpha_i)(\ul\phi^{i,p})\}\subset\{\ul\phi^{i,p-1}\}+V_{<\alpha_i}$ for all $(i,p)\in[1,\nu]\times\Z$.
\end{enumerate}
\item[Output:] A matrix $M\in\Ctt^{\mu\times\mu}$ such that $\ul\phi M$ is a good basis of $L$.
\item Compute a minimal Gr\"obner basis
\[
M:=\GB(M)\in\Cth^{\mu\times\mu}
\]
compatible with the ordering $>$ on $\{\theta^k\ul\phi^{i,p}\mid(k,i,p)\in\Z\times[1,\nu]\times\Z\}$ defined by
\[
(k,i,p)>(l,j,q):\Leftrightarrow k>l\lor(k=l\land(i>j\lor(i=j\land p>q)))
\]
for all $(k,i,p),(l,j,q)\in\Z\times[1,\nu]\times\Z$.
\item Set $(M^{k,i})_{(k,i)\in\Z\times[1,\nu]}:=M$ where
\[
\{\lexp(M^{k,i})\}=\{(k,i,n-k)\}
\]
for all $(k,i)\in\Z\times[1,\nu]$.
\item\label{7} Compute $(A^{k,i}_{s,l,j})_{(s,j,l)\in\Z\times[1,\nu]\times\Z}$ and 
\begin{align*}
\Phi^{k,i}_{s,l,j}&:=\theta(B-(\alpha_i+k)+\theta\dth)M^{k,i}-M^{1+k,i}A^{k,i}_{0,1+k,i}\\
&-\sum_{\substack{(l',j')<(1+k,i)\\(s+l,j,n-l)<(l',j',n-l')}}M^{l',j'}A^{k,i}_{0,l',j'}-\theta^{s}M^{l,j}A^{k,i}_{s,l,j}
\end{align*}
such that $\lexp(\Phi^{k,i}_{s,l,j})<(s+l,j,n-l)$ for all $(k,i)\in\Z\times[1,\nu]$ for decreasing $(s+l,j,n-l)$ until $\Phi^{k,i}_{s,l,j}=0$ or $A^{k,i}_{s,l,j}\ne0$ and $s\ge1$.
\item If $\Phi^{k,i}_{s,l,j}=0$ then return $M:=(M^{k,i})_{(k,i)\in\Z\times[1,\nu]}$.
\item Choose $(k,i)\in\Z\times[1,\nu]$ and $(s,j,l)\in\Z\times[1,\nu]\times\Z$ with $A^{k,i}_{s,l,j}\ne0$ and $s\ge1$ such that $(s+l,j,n-l)$ is maximal.
\item Set $c^{k,i}_{s,l,j}:=(1+k+\alpha_i-s-l-\alpha_j)^{-1}$ and
\[
M^{k,i}:=M^{k,i}+c^{k,i}_{s,l,j}\theta^{s-1}M^{l,j}A^{k,i}_{s,l,j}
\]
\item Go to (\ref{7}).
\end{enumerate}
\end{alg}

\begin{lem}\label{32}
Algorithm \ref{9} terminates and is correct.
\end{lem}
\begin{pf}
By Lemma \ref{3},
\[
V_{\alpha_j+q}G=\bigoplus_{(p,i)\le(q,j)}\theta^p\langle\ul\phi^i\rangle\C
\]
for all $(q,j)\in\Z\times[1,\nu]$.
Let $0\ne\ol m\in\langle M\rangle\Cth$ with $\lexp(\ol m)=(k,i,p)$.
Then
\[
\tau^k\ul\phi\ol m\in\langle\ul\phi^{i,q}\mid q\le p\rangle\C+V_{<\alpha_i}G
\]
and, since $\ul\phi\ol m\in\langle\ul\phi M\rangle\Cth=L$ and $L^{n-k}\gr^V_{\alpha_i}G=\langle\ul\phi^{i,q}\mid q\ge n-k\rangle\C$,
\[
\tau^k\ul\phi\ol m\in\langle\ul\phi^{i,q}\mid q\ge n-k\rangle\C+V_{<\alpha_i}G.
\]
In particular, $p\ge n-k$.
Moreover, $\tau^k\ul\phi\lead(\ol m)\in\gr_L^p\gr^V_{\alpha_i}G$ and hence, by Lemma \ref{6},
\[
\theta^{n-p-k}\ul\phi\lead(\ol m)\in\gr_L\gr^VL=\langle\ul\phi\lead(M)\rangle\Cth.
\]
In particular, if $p>n-k$ then $\lead(\ol m)\in\theta\langle\lead(M)\rangle\Cth$.
Since $M$ is a minimal Gr\"obner basis, this implies that
\[
\{\lexp(M^{k,i})\}=\{(k,i,n-k)\},\quad M^{k,i}\equiv\lead(M^{k,i})\mod\text{terms}<(k,i)
\]
for all $(k,i)\in\Z\times[1,\nu]$.
In particular, $\{\ul\phi M^{k,i}\}\subset V_{k+\alpha_i}$ for all $(k,i)\in\Z\times[1,\nu]$.

By Lemma \ref{21}, $t\circ\ul\phi=\theta(-\tau\dta)\circ\ul\phi=\ul\phi\circ\theta(B+\theta\dth)$.
Since 
\[
\{\theta(\theta\dth-(k+\alpha_i))(\theta^k\ul\phi^{i,n-k})\}\subset\{\theta^{1+k}\ul\phi^{i,n-(1+k)}\}\mod V_{<1+k+\alpha_i},
\]
there is a matrix $A^{i,k}_{0,1+k,i}$ such that
\[
\theta(B-(k+\alpha_i)+\theta\dth)M^{k,i}\equiv M^{1+k,i}A^{i,k}_{0,1+k,i}\mod\text{terms}<(1+k,i)
\]
and hence there are matrices $A^{i,k}_{s,l,j}$ such that
\[
\theta(B-(k+\alpha_i)+\theta\dth)M^{k,i}-M^{1+k,i}A^{i,k}_{0,1+k,i}=
\sum_{(s'+l',j')<(1+k,i)}\theta^{s'}M^{l',j'}A^{k,i}_{s',l',j'}
\]
for all $(k,i)\in\Z\times[1,\nu]$.
Choose $(k,i)\in\Z\times[1,\nu]$ and $(s,j,l)\in\Z\times[1,\nu]\times\Z$ such that $(s+l,j,n-l)$ is maximal with $A^{k,i}_{s,l,j}\ne0$ and $s\ge1$.
In particular, $(s+l,j)<(1+k,i)$ and hence $1+k+\alpha_i-s-l-\alpha_j>0$ and $c^{k,i}_{s,l,j}>0$ is defined.
Moreover, since
\[
\{(\theta\dth-(\alpha_j+l))(\theta^l\ul\phi^{j,n-l})\}\subset\{\theta^l\ul\phi^{j,n-(1+l)}\}\mod V_{<\alpha_j+l},
\]
\begin{align*}
\Phi^{k,i}_{s,l,j}&=\theta(B-(k+\alpha_i)+\theta\dth)\bigl(M^{k,i}-c^{k,i}_{s,l,j}\theta^{s-1}M^{l,j}A^{k,i}_{s,l,j}\bigr)\\
&-M^{1+k,i}A^{i,k}_{0,1+k,i}-\sum_{\substack{(l',j')<(1+k,i)\\(s+l,j,n-l)<(l',j',n-l')}}M^{l',j'}A^{k,i}_{0,l',j'}\\
&\equiv\theta^sM^{l,j}A^{k,i}_{s,l,j}
+c^{k,i}_{s,l,j}\theta(\theta\dth-k-\alpha_i)\theta^{s-1}M^{l,j}A^{k,i}_{s,l,j}\\
&\equiv\theta^sM^{l,j}A^{k,i}_{s,l,j}
+c^{k,i}_{s,l,j}\theta^s(\theta\dth+s-1-k-\alpha_i)M^{l,j}A^{k,i}_{s,l,j}\\
&\equiv\theta^sM^{l,j}A^{k,i}_{s,l,j}
+c^{k,i}_{s,l,j}\theta^s(s+l+\alpha_j-1-k-\alpha_i)M^{l,j}A^{k,i}_{s,l,j}\\
&\equiv0\mod\text{ terms}<(s+l,j,n-l)
\end{align*}
and hence $(s+l,j,n-l)$ is strictly decreasing until $\Phi^{k,i}_{s,l,j}=0$.
Then the algorithm terminates and $\ul\phi M=(\ul\phi M^{k,i})_{(k,i)\in\Z\times[1,\nu]}$ is a $\Cth$--basis of $L$ with 
\[
t(\ul\phi M^{k,i})=\ul\phi M^{1+k,i}A^{i,k}_{0,1+k,i}+
\sum_{(l',j')<(1+k,i)}\ul\phi M^{l',j'}A^{k,i}_{0,l',j'}+\theta(k+\alpha_i)\ul\phi M^{k,i}
\]
and $\{\ul\phi M^{k,i}\}\subset V_{k+\alpha_i}$ for all $(k,i)\in\Z\times[1,\nu]$.
Hence, $\ul\phi M$ is a good basis of $L$.
\end{pf}

The following algorithm combines Algorithms \ref{8}, \ref{1}, \ref{22}, \ref{4}, and \ref{9} to compute a good basis of $G_0$.

\begin{alg}\label{25}\
\begin{enumerate}
\item[Input:]
A cohomologically tame polynomial $f\in\CX$.
\item[Output:]
\begin{enumerate}
\item A vector $\ul\phi\in\CXth^\mu$ such that $[\ul\phi]$ is a good basis of $G_0$.
\item The matrix $A=A^{[\ul\phi]}\in\Cth^{\mu\times\mu}$ of $t$ with respect to $[\ul\phi]$.
\end{enumerate}
\item Set $k:=\deg(f)$.
\item\label{24} Compute $\ul\phi\in\CXth^\mu$ and $A\in\Cth^{\mu\times\mu}$ by Algorithm \ref{8}.
\item Compute $U\in\Cth^{\mu\times\mu}$ and $B\in\Cta^{\mu\times\mu}$ by Algorithm \ref{1}.
\item Set $\ul\phi:=\ul\phi U$, $B:=\inv U(B-\tau\dta)U\in\Cta^\mu$, and $M:=\inv U\in\Ctt^\mu$.
\item Compute $\sigma$ by Algorithm \ref{22}.
\item If $\frac{1}{\mu}\sum\sigma>\frac{n+1}{2}$ then set $k:=k+1$ and go to (\ref{24}).
\item Compute $U=(U^{i,p})_{(i,p)\in[1,\nu]\times\Z}\in\GL_\mu(\C)$ by Algorithm \ref{4}.
\item Set $\ul\phi:=(\ul\phi U^{i,p})_{(i,p)\in[1,\nu]\times\Z}$, $B:=\inv U(B-\tau\dta)U\in\Cta^{\mu\times\mu}$, and $M:=\inv UM$.
\item Compute $M\in\Ctt^{\mu\times\mu}$ by Algorithm \ref{9}.
\item Set $\ul\phi:=\ul\phi M$ and $A:=\inv M\theta(B-\tau\dta)M\in\Cth^{\mu\times\mu}$.
\item Return $\ul\phi$ and $A$.
\end{enumerate}
\end{alg}

\begin{prop}
Algorithm \ref{25} terminates and is correct.
\end{prop}
\begin{pf}
Let $L_k\subset G_0$ be computed by Algorithm \ref{8}.
Then $L_k=G_0$ for $k\gg0$ and $k$ is strictly increasing while $\frac{1}{\mu}\sum\sigma>\frac{n+1}{2}$.
By Lemma \ref{27} and Theorem \ref{26}, $L=G_0$ if and only if $\frac{1}{\mu}\sum\sigma=\frac{1}{\mu}\sum\spec(L)=\frac{1}{\mu}\sum\spec(G_0)=\frac{n+1}{2}$.
This implies that $L_k=G_0$ after finitely many steps.
By Theorem \ref{23}, $L:=L_k=G_0$ is a good lattice as required by algorithms \ref{4} and \ref{9}.
Hence, the algorithm terminates and is correct.
\end{pf}

\begin{rem}
In the local situation, one can replace the algorithms \cite[7.4--5]{Sch02c} by the algorithms \ref{4} and \ref{9} to avoid the linear algebra computation \cite[7.4]{Sch02c}.
This modified algorithm is implemented in the {\sc Singular} \cite{GPS04} library {\tt gmssing.lib} \cite{Sch04c}.
\end{rem}

\section{Examples}

Algorithm \ref{25} is implemented in the {\sc Singular} \cite{GPS04} library {\tt gmspoly.lib} \cite{Sch04d}.
Using this implementation, we compute a good basis $\ul\phi$ of $G_0$ for several examples.
By Lemma \ref{37}, the diagonal of $A^\ul\phi_1$ determines the spectrum of $f$.
Using Proposition \ref{34}, we read off the monodromy $T_\infty$ around the discriminant of $f$ from $A^\ul\phi$.
First, we compute two convenient and Newton non--degenerate examples \cite{Dou99}.

\begin{exmp}
Let $f=x^2+y^2+x^2y^2$. Then {\sc Singular} computes
\[
\ul\phi=\Bigl(1,xy,y,x,x^2+\frac{1}{2}\Bigr)
\]
and 
\[
A^{[\ul\phi]}=\begin{pmatrix}
-\frac{1}{2}&0&0&0&\frac{1}{4}\\
0&-1&0&0&0\\
0&0&-1&0&0\\
0&0&0&-1&0\\
1&0&0&0&-\frac{1}{2}
\end{pmatrix}+
\theta\begin{pmatrix}
\frac{1}{2}&0&0&0&0\\
0&1&0&0&0\\
0&0&1&0&0\\
0&0&0&1&0\\
0&0&0&0&\frac{3}{2}
\end{pmatrix}.
\]
The monodromy $T_\infty$ has a $2\times 2$ Jordan block with eigenvalue $-1$. 
\end{exmp}

\begin{exmp}
Let $f=x+y+z+x^2y^2z^2$. Then {\sc Singular} computes
\[
\ul\phi=\Bigl(1,\theta^2x-3\theta x^2+x^3,\frac{5}{2}x,10\theta^2x^2-\frac{25}{2}\theta x^3+\frac{5}{2}x^4,-\frac{25}{4}\theta x+\frac{25}{4}x^2\Bigr)
\]
and 
\[
A^{[\ul\phi]}=\begin{pmatrix}
0&0&0&-\frac{25}{8}&0\\
0&0&0&0&\frac{125}{8}\\
1&0&0&0&0\\
0&1&0&0&0\\
0&0&1&0&0
\end{pmatrix}+
\theta\begin{pmatrix}
\frac{1}{2}&0&0&0&0\\
0&1&0&0&0\\
0&0&\frac{3}{2}&0&0\\
0&0&0&2&0\\
0&0&0&0&\frac{5}{2}
\end{pmatrix}.
\]
The monodromy $T_\infty$ has a $2\times 2$ Jordan block with eigenvalue $1$ and a $3\times 3$ Jordan block with eigenvalue $-1$. 
\end{exmp}

Finally, we compute a non--convenient and Newton degenerate but tame \cite[3]{Bro88} example.

\begin{exmp}
Let $f=x(x^2+y^3)^2+x$. Then {\sc Singular} computes
\begin{multline*}
\ul\phi=\Bigl(
1,
623645y,
\frac{8645}{24}x,
-2470\bigl(y^3+x^2\bigr),\\
-11339\bigl(y^4+x^2y\bigr),
475\bigl(2\theta y^3-5\theta x^2+6x^3\bigr),
y^2,
3\theta^2y^2+4y^5,\\
6670\bigl(\theta y^4-10\theta x^2y+6x^3y\bigr),
8\theta^2y^3-20\theta^2x^2-15\theta x^3+18x^4+3,\\
-4365515\bigl(35\theta^2y^4-350\theta^2x^2y-300\theta x^3y+180x^4y+24y\bigr),
\frac{623645}{6}xy,\\
-8645\bigl(\theta x+2y^3-4x^2\bigr),
-124729\bigl(5\theta xy+y^4-5x^2y\bigr)
\Bigr)
\end{multline*}
and $A^{[\ul\phi]}=A^{[\ul\phi]}_0+\theta A^{[\ul\phi]}_1$ where
{\tiny
\[
\setcounter{MaxMatrixCols}{20}
A^{[\ul\phi]}_0=\begin{pmatrix}
0&0&0&0&0&-380&0&0&0&0&0&0&0&0\\
0&0&0&0&0&0&0&0&-\frac{32}{4675}&0&0&0&0&0\\
\frac{96}{43225}&0&0&0&0&0&0&0&0&\frac{-288}{216125}&0&0&0&0\\
0&0&-\frac{7}{180}&0&0&0&0&0&0&0&0&0&0&0\\
0&0&0&0&0&0&0&0&0&0&0&-\frac{11}{9}&0&0\\
0&0&0&-\frac{52}{75}&0&0&0&0&0&0&0&0&\frac{728}{75}&0\\
0&0&0&0&0&0&0&0&0&0&0&0&0&0\\
0&0&0&0&0&0&0&0&0&0&0&0&0&0\\
0&0&0&0&-\frac{17}{75}&0&0&0&0&0&0&0&0&\frac{187}{15}\\
0&0&0&0&0&\frac{380}{3}&0&0&0&0&0&0&0&0\\
0&0&0&0&0&0&0&0&-\frac{4}{98175}&0&0&0&0&0\\
0&\frac{24}{5}&0&0&0&0&0&0&0&0&\frac{2016}{5}&0&0&0\\
0&0&\frac{1}{180}&0&0&0&0&0&0&0&0&0&0&0\\
0&0&0&0&0&0&0&0&0&0&0&\frac{1}{9}&0&0
\end{pmatrix}
\]}
and
{\tiny
\[
A^{[\ul\phi]}_1=\begin{pmatrix}
\frac{1}{3}&0&0&0&0&0&0&0&0&0&0&0&0&0\\
0&\frac{7}{15}&0&0&0&0&0&0&0&0&0&0&0&0\\
0&0&\frac{2}{3}&0&0&0&0&0&0&0&0&0&0&0\\
0&0&0&\frac{11}{15}&0&0&0&0&0&0&0&0&0&0\\
0&0&0&0&\frac{13}{15}&0&0&0&0&0&0&0&0&0\\
0&0&0&0&0&\frac{14}{15}&0&0&0&0&0&0&0&0\\
0&0&0&0&0&0&1&0&0&0&0&0&0&0\\
0&0&0&0&0&0&0&1&0&0&0&0&0&0\\
0&0&0&0&0&0&0&0&\frac{16}{15}&0&0&0&0&0\\
0&0&0&0&0&0&0&0&0&\frac{17}{15}&0&0&0&0\\
0&0&0&0&0&0&0&0&0&0&\frac{19}{15}&0&0&0\\
0&0&0&0&0&0&0&0&0&0&0&\frac{4}{3}&0&0\\
0&0&0&0&0&0&0&0&0&0&0&0&\frac{23}{15}&0\\
0&0&0&0&0&0&0&0&0&0&0&0&0&\frac{5}{3}
\end{pmatrix}.
\]}
The monodromy $T_\infty$ is unipotent with eigenvalues
\begin{gather*}
\ee^{-2\pi\ii\frac{1}{3}},\ee^{-2\pi\ii\frac{7}{15}},\ee^{-2\pi\ii\frac{2}{3}},\ee^{-2\pi\ii\frac{11}{15}},\ee^{-2\pi\ii\frac{13}{15}},\ee^{-2\pi\ii\frac{14}{15}},1,\\1,\ee^{-2\pi\ii\frac{16}{15}},\ee^{-2\pi\ii\frac{17}{15}},\ee^{-2\pi\ii\frac{19}{15}},\ee^{-2\pi\ii\frac{4}{3}},\ee^{-2\pi\ii\frac{23}{15}},\ee^{-2\pi\ii\frac{5}{3}}.
\end{gather*}
\end{exmp}

\bibliographystyle{elsart-num}
\bibliography{gbtp}

\begin{thebibliography}{10}
\expandafter\ifx\csname url\endcsname\relax
  \def\url#1{\texttt{#1}}\fi
\expandafter\ifx\csname urlprefix\endcsname\relax\def\urlprefix{URL }\fi

\bibitem{Sab98b}
C.~Sabbah, Hypergeometric periods for a tame polynomial, arXiv.org
  math.AG/9805077.

\bibitem{Bri70}
E.~Brieskorn, {Die Monodromie der isolierten Singularit\"aten von
  Hyperfl\"achen}, Manuscr. Math. 2 (1970) 103--161.

\bibitem{Seb70}
M.~Sebastiani, Preuve d'une conjecture de {Brieskorn}, Manuscr. Math. 2 (1970)
  301--308.

\bibitem{Ste76}
J.~Steenbrink, Mixed {Hodge} structure on the vanishing cohomology, in: Real
  and complex singularities, Nordic summer school, Oslo, 1976, pp. 525--562.

\bibitem{Pha79}
F.~Pham, Singularit{\'e}s des syst{\`e}mes de {Gauss}--{Manin}, Vol.~2 of
  Progr. in Math., Birkh{\"a}user, 1979.

\bibitem{Var82a}
A.~Varchenko, Asymptotic {Hodge} structure in the vanishing cohomology, Math.
  USSR Izv. 18~(3) (1982) 496--512.

\bibitem{Sai89}
M.~Saito, On the structure of {Brieskorn} lattices, Ann. Inst. Fourier Grenoble
  39 (1989) 27--72.

\bibitem{Sch02c}
M.~Schulze, The differential structure of the {Brieskorn} lattice, in:
  A.~Cohen, et~al. (Eds.), Mathematical Software --- ICMS 2002, World Sci.,
  2002, pp. 136--146.

\bibitem{Sch04a}
M.~Schulze, A normal form algorithm for the {Brieskorn} lattice, to appear in
  J. Symb. Comp. .

\bibitem{Dou99}
A.~Douai, Tr\`es bonnes bases du r\'eseau de {Brieskorn} d'un polyn\^ome
  mod\'er\'e, Bull. Soc. Math. France 127 (1999) 255--287.

\bibitem{KV85}
A.~Khovanskii, A.~Varchenko, Asymptotics of integrals over vanishing cycles and
  the {Newton} polyhedron, Sov. Math. Docl. 32 (1985) 122--127.

\bibitem{Dou93}
A.~Douai, {\'E}quations aux diff\'erences finies, int\'egrales de fonctions
  multiformes et poly\`edres de {Newton}, Comp. math. 87 (1993) 311--355.

\bibitem{BGM89}
J.~Brian{\c{c}}on, M.~Granger, P.~Maisonobe, M.~Miniconi, Algorithme de calcul
  du polyn\^ome de {Bernstein}, Ann. Inst. Fourier 3 (1989) 553--609.

\bibitem{Bor87}
A.~Borel, et~al. (Eds.), Algebraic {D}--modules, 2nd Edition, Vol.~2 of Persp.
  in Math., Acad. Press, 1987.

\bibitem{Eis96}
D.~Eisenbud, Commutative Algebra with a View toward Algebraic Geometry, Vol.
  150 of Grad. Texts in Math., Springer, 1996.

\bibitem{Sab98a}
C.~Sabbah, Monodromy at infinity and {Fourier} transform, Publ. RIMS, Kyoto
  Univ. 33 (1998) 643--685.

\bibitem{Sab02}
C.~Sabbah, D\'eformations isomonodromiques et vari\'et\'es de {Frobenius}, EDP
  Sciences, 2002.

\bibitem{GP02}
G.~Greuel, G.~Pfister, A {\sc Singular} Introduction to Commutative Algebra,
  Springer, 2002.

\bibitem{Sab87}
C.~Sabbah, $\mathcal{D}$-modules et cycles {\'e}vanescents, in: J.-M. Aroca,
  T.~Sanchez-Giralda, J.-L. Vicente (Eds.), Deuxi{\`e}me conf{\'e}rence de {La
  Rabida}, G{\'e}om{\'e}trie algebrique et applications III, Vol.~24 of Travaux
  en cours, Hermann, Paris, 1987, pp. 53--98.

\bibitem{GPS04}
G.-M. Greuel, G.~Pfister, H.~Sch{\"o}nemann, {\sc Singular} 2.0.5, {A Computer
  Algebra System for Polynomial Computations}, Centre for Computer Algebra,
  University of Kaiserslautern, {\tt http://www.singular.uni-kl.de} (2004).

\bibitem{Sch04c}
M.~Schulze, {\tt gmssing.lib}, {\sc Singular} library, Centre for Computer
  Algebra, University of Kaiserslautern, {\tt http://www.singular.uni-kl.de}
  (2004).

\bibitem{Sch04d}
M.~Schulze, {\tt gmspoly.lib}, {\sc Singular} library, Centre for Computer
  Algebra, University of Kaiserslautern, {\tt http://www.singular.uni-kl.de}
  (2004).

\bibitem{Bro88}
S.~Broughton, {Milnor} numbers and the topology of polynomial hypersurfaces,
  Inv. Math. 92 (1988) 217--241.

\end{thebibliography}

\end{document}